\documentclass[12pt]{amsart}

\usepackage{amssymb,latexsym}

\usepackage{enumerate}

\makeatletter

\@namedef{subjclassname@2010}{%

  \textup{2010} Mathematics Subject Classification}

\makeatother

\makeatletter

\@namedef{subjclassname@2010}{%

  \textup{2010} Mathematics Subject Classification}

\makeatother

\newtheorem{theorem}{Theorem}[section]

\newtheorem{corollary}[theorem]{Corollary}

\newtheorem{lemma}[theorem]{Lemma}
\newtheorem{proposition}[theorem]{Proposition}

\theoremstyle{definition}

\numberwithin{equation}{section}

\newcommand{\subgp}[1]{\langle{#1}\rangle}

\begin{document}


\baselineskip=17pt



\title{Two Inverse results  }
\author{ {Y. O. Hamidoune}\\
UPMC, Univ Paris 06,\\ 4 Place Jussieu,
75005 Paris, France.}
\email{hamidoune@math.jussieu.fr}

\maketitle

\begin{abstract}
Let $ A$ be a subset of group $G_0$
 with $|{A^{-1}A}|\le 2|A|-2.$ We show that there are an element $a\in A$  and a non-null proper subgroup $H$ of $G$ such that one of the following holds:
\begin{itemize}
  \item  $x^{-1}Hy \subset A^{-1}A,$ for all $(x,y)\in A^2\setminus (Ha)^2,$
  \item  $xHy^{-1} \subset AA^{-1},$ for all $(x,y)\in A^2\setminus (aH)^2.$ 
\end{itemize}
where  $G$ is the subgroup generated by ${A^{-1}A}.$

Assuming that $A^{-1}A\neq G$ and that $ |A^{-1}A|< \frac{5|A|}3,$
we show that there are a normal subgroup $K$ of  $G$ and a subgroup  $H$ with $K\subset H\subset A^{-1}A $ and  $2|K|\ge |H|$ such that
$$A^{-1}AK=KA^{-1}A=A^{-1}A\  \text{and}\   6|K|\ge |A^{-1}A|=3|H|.$$
\end{abstract}

\section{Introduction}

Let  $ A,B$ be
subsets of a group $ G $. The {\em Minkowski product} of $ A$ with $B$ is defined as
$$AB=\{xy \ : \ x\in A\  \mbox{and}\ y\in
  B\}.$$

 The Cauchy-Davenport Theorem \cite{cauchy,davenport} states that  $|AB|\ge |A|+|B|-1,$ if $AB$ is a proper subset of $G$ and if $|G|$ is
 a prime.  Kneser's Theorem \cite{knesrcomp} states that  $AB$ is a periodic set if $|AB|\le |A|+|B|-2$ and if  $G$ is abelian. Diderrich \cite{diderrich} obtained the same conclusion assuming only that the elements of $B$ commute. As mentioned in \cite{hfour}, the last result follows from Kneser's Theorem. In \cite{olsonsdif}, Olson constructed subsets $A$ and $B$ of some non-abelian group with $|AB|\le |A|+|B|-2$ such that for every  non-null group $H,$

  $AB\neq AHB$,  $AB\neq HAB$ and $AB\neq ABH.$

  Also, the  cases $B=A$ and $B=A^{-1}$ received also some attention.
 In \cite{freiman},
 Freiman obtained an inverse result describing  $A$ if $|A^2|<1.6|A|.$
 A transparent  exposition of Freiman results is contained in Husbands dissertation \cite{husb}. In the last dissertation,
one may find a lemma due to Freiman,  allowing to recover inverse results for
$A^2$ from inverse results $A A^{-1}.$ Recently, Tao \cite{t2},  gave a short argument proving an inverse theorem for
$|A^{-1} A|< \frac{1+\sqrt{5}}{2}|A|.$ Tao adds in \cite{t2}:

  "One should be able to get a bit more structural information on ${A}$ than is given by the above conclusion, and I doubt the golden ratio is sharp either  (the correct threshold should be {2}, in analogy with the commutative Kneser theorem; after that, the conclusion will fail, as can be seen by taking ${A} $ to be a long geometric progression). Readers here are welcome to look for improvements to these results, of course."

  Our first result the following Kneser's type theorem:

 Let $ A$ be a subset of group $G_0$  with $|A^{-1}A|\le 2|A|-2$ and $A^{-1}A\neq G,$
where  $G$ is the subgroup generated by ${A^{-1}A}.$
We show that  there are an element $a\in A$  and a non-null subgroup $H$ of $G$ such that one of the following holds:
\begin{itemize}
  \item  For all $(x,y)\in A^2\setminus (Ha)^2,$  $x^{-1}Hy  \subset A^{-1}A,$ and $|A^{-1}A|>|A^{-1}H|+|HA|-2|H|.$
  \item  For all $(x,y)\in A^2\setminus (aH)^2,$  $xHy^{-1}  \subset AA^{-1},$ and  $|AA^{-1}|>|AH|+|AH|-2|H|.$
\end{itemize}

  Improving the bounds in literature, we prove the following:

Let $ A$ be a subset of group $G_0$ and let $G$ be the subgroup generated by
${A^{-1}A}.$  If $|A^{-1}A|<\min (|G|, \frac{5|A|}3),$
then there are a normal subgroup
 $K$ of  $G$ and a subgroup  $H$ with $K\subset H\subset A^{-1}A $ and  $2|K|\ge |H|$ such that
$$A^{-1}AK=KA^{-1}A=A^{-1}A \  \text{and}\   6|K|\ge|A^{-1}A|=3|H|.$$

In \cite{freiman}, Freiman proved that $A^2=xA\cdot A^{-1},$ for some $x\in A^2,$ if $|A^2|< \frac{1+\sqrt{5}}{2}|A|.$
Freiman's result may be combined with our last result to deduce  an inverse theorem for $|A^2|<\frac{1+\sqrt{5}}{2}|A|,$
improving slightly  the inverse result for $|A^2|<1.6|A|$ mentioned above.

Let $G$ be a group,  $H$ be a subgroup and $a$ be an element of $G\setminus H$
such that $H\cup \{a\}$ generates a subgroup with order $>2|H|.$
Put $E=H\cup Ha.$ Clearly, $E^{-1}E$ is the union of three $H$--cosets,
if $a\in N(H),$ where $N(H)$ is the normalizer of $H.$

If $a\notin N(H),$  and if $H$ has a subgroup $K$ with $2|K|=|H|,$ which is normal in $G.$
Clearly $E^{-1}E$ is not a union of $5$ cosets of the same subgroup, if $|H|$ is large.  This example  shows that
the bounds in Theorem \ref{periodic} are somehow tight.

Let $S$ be a generating subset of a group $G$ containing  $1$
let $k$ be an integer such that  $|S|\ge k+1.$
As a first step in the proof of Theorem \ref{periodic}, Proposition \ref{lesstwothird} states
that if for some proper subset $T\subset G,$  $|TS|<\min (|G|,|T|+(1-\frac{1}{k})|S|)$ holds, then there is a finite  subgroup $H$
 such that
$|HS| \le (k-1)|H|$ or $ |SH|\le (k-1)|H|.$ In particular, $A$ is covered by few cosets if $A^{-1}A$ (resp. $A^2$ ) has a small cardinality.
Quite likely, this conclusion may be used to describe subsets $A$ with $ |A^2|< \frac{5|A|}3.$ It could be also useful in the description of the subsets $A$ with $ |A^{-1}A|< \frac{9|A|}5.$

We use the isoperimetric approach, which could be appropriate in the  investigation of some inverse problems.
Our main tool is a result  proved by the author in \cite{hejc}, stating that for some $T\in \{S,S^{-1}\},$ the objective function
$|XT|-|X|,$ where $1\leq |XT|<|G|,$ achieves its minimum value on a subgroup. In order to make the present work self-contained, we include a proof of this result.

\section{Preliminaries}

Let $A$ be a subset of a group  $G$ and let $H$ be a subgroup. By a $H$-{\em right-component} of $A,$
we shall mean a non-empty trance on $A$ of some right $H$-coset. It is thus, a set of the form $A\cap (H+x),$
where $x\in A.$

  Recall a well known fact:

\begin{lemma}(folklore)
Let $a,b$ be elements of a group $G$ and let $H$ be a finite subgroup of $G.$
Let
 $A,B$  be  subsets of $G$
 such that  $A\subset aH$ and $B\subset Hb.$ If  $|A|+|B|>|H|,$
 then $AB=aHb$.
\label{prehistorical}
 \end{lemma}


\begin{lemma} \label{Obs0}
Let $a$ be an element of a group $G$ and let $A$ be a finite   subset of $G.$

 Then $a^{-1}\subgp{A^{-1}A}a=\subgp{(Aa)^{-1}Aa}.$ Moreover $K$ is a normal subgroup of $\subgp{A^{-1}A}$ if and only if
 $a^{-1}Ka$ is a normal subgroup of $\subgp{(Aa)^{-1}Aa}.$
\end{lemma}

\begin{lemma} \label{Obs}

Let $x,y$ be elements of a group $G$ and  let $H$ be a subgroup of $G.$
For every $c\in xH\cap Hy,$

\begin{equation}\label{obs}
 xH\cap Hy=((xH x^{-1})\cap H)c=c((y^{-1}Hy)\cap H).
 \end{equation}


 In particular, if for some $a\in G,$
$Ha\cap a^{-1}H\neq \emptyset,$ and if $H\cap (a^{-1}Ha)$
is a normal subgroup of $H,$ then $H\cap (a^{-1}Ha)$
is a normal subgroup of $\subgp{H\cup \{a\}}.$
\end{lemma}

\begin{proof}

Since $xH=cH$ and $Hy=Hc,$ we have
 $$xH\cap Hy=cH\cap Hc=((cHc^{-1})\cap H)c=((xHx^{-1})\cap H)c.$$
Similarly, $$xH\cap Hy=cH\cap Hc=c((c^{-1}Hc)\cap H)=c((y^{-1}Hy)\cap H).$$

Assume now that $Ha\cap a^{-1}H\neq \emptyset$  and put $K=H\cap (a^{-1}Ha).$

Choose an element $c\in a^{-1}H\cap Ha$ and put  $c=ea,$ for some $e\in H.$ By (\ref{obs}),
$Kc=cK.$  Thus
$eKa=ea K,$  and hence $Ka=aK.$ Therefore,  $K$ is a normal subgroup of $\subgp{H\cup \{a\}}.$\end{proof}

Lemmas \ref{Obs0} and \ref{Obs} are easy exercises.

We shall give succinct presentation of the isoperimetric approach making the present work self-contained.

Let $G$ be a group and let $S$ be a generating finite subset with $1\in S.$ For a subset $X\subset G,$ we set $\nabla (X)=XS\setminus X.$
 Put $${\mathcal F}(S)=\{X\subset G : XS\ \text{is a proper finite  subset of} \ G.\} $$
The {\em connectivity} of $S$ is  $$\kappa (S)=\min \{|XS|-|X| : X\in {\mathcal F} (S)\},$$
 where $\min (\emptyset)=|G|.$ A member of ${\mathcal F}(S)$ achieving this minimum will be called a fragment of $S$. A fragment with minimal cardinality
will be called an atom.

We note that for $1\in S,$ one has $\subgp{Sa^{-1}}=\subgp{a^{-1}S}=\subgp{S}$
and that $$\kappa (Sa^{-1})=\kappa ({S})=\kappa ({a^{-1}S}),$$ for every $a\in S.$
The next lemma is a basic one in order to deal with finite groups where a  switch from $S$ to $S^{-1}$ is necessary.

A subset $S$ will be called {\em faithful} if $|A|\le |\nabla (A)|.$
In particular, all finite subsets are faithful if $G$ is infinite.

\begin{lemma} \label{finiteg}{Let $S$ be a generating  subset of a finite group $G$ with $1\in S.$
Then \begin{equation}
\kappa (S)= \kappa (S^{-1}).\label{kk-}\end{equation}  Moreover,
\begin{itemize}
 \item A subset $X\subset G$ {is a } {fragment } of $S$ {if and only if } $\nabla(X)$ is a fragment of $S^{-1}.$
 \item If $S$ is non-faithful, then $G$ is a non-abelian group and
 $S^{-1}$ is faithful.
 \end{itemize}
   }\end{lemma}
\begin{proof}
Notice that $\nabla (X)S^{-1}\cap X=\emptyset,$ and hence $\nabla (X)S^{-1} \subset G\setminus X.$
In particular, $\nabla (X)\in  {\mathcal F}(S^{-1})$ if and only if $X\in  {\mathcal F}(S^{-1}).$

In order to prove (\ref{kk-}), we may  assume that $ {\mathcal F}(S)\neq \emptyset,$ otherwise
 (\ref{kk-}) holds by convention. Take a fragment $X.$ We have
\begin{align*} \kappa (S^{-1})&\le
|\nabla (X)S^{-1}|-|\nabla (X)|\\ &\le |G|-|X|-|\nabla (X)|=\kappa (S).\end{align*}
The reverse inequality follows by duality. Thus, the inequalities chain consists of  equalities, and thus
$\nabla(X)$ is a fragment of $S^{-1}.$

Assume that $S$ is non-faithful and let $K$ be an atom of $S^{-1}.$
 By the first part of the Lemma,
$\nabla (H)$ is a fragment of $S^{-1}.$
Thus, $|H|> |\nabla (H)|\ge |K|.$
By  (\ref{kk-}),
$$|K|+\kappa (S^{-1})+|\nabla^- (K)|=|G|=|H|+\kappa (S)+|\nabla (H)|,$$
we have
 $ |K|\le |\nabla (H)|<|\nabla^- (K)|.$
In particular, $S^{-1}$ is faithful.

In the abelian case, $H^{-1}$ is an  atom of $S^{-1}$. Hence, $S$ and $S^{-1}$ are
are both faithfull.\end{proof}

\begin{theorem}
Let $S$ be a generating finite  subset of a group $G$ with $1\in S.$

Let  $A$ be an atom of $S$ and let $F$ be a fragment of
    $S.$ If $|A|\le | \nabla (F)|$ and $A\cap F\neq \emptyset,$ then $A\subset F.$
    In particular, distinct atoms have an empty intersection, if $S$ is faithful.
\label{lem1977} \end{theorem}

\begin{proof}
Assume that $A\cap F\neq \emptyset.$

{\begin{center}
\begin{tabular}{|c||c|c|c|c|}
\hline
$\cap $&  $F$ & $FS\setminus F $ &  $\nabla (F)$ \\
\hline
\hline
$A$&  $R_{11}$ & $R_{12}$ & $R_{13}$ \\
\hline
$AS\setminus A $& $R_{21}$& $R_{22}$ & $R_{23}$ \\
\hline
$\nabla (A)$ & $R_{31}$  & $R_{32}$ & $R_{33}$ \\
\hline
\end{tabular}
\end{center}}

By the definition of $\kappa,$
\begin{align*}  |R_{21}|+|R_{22}|+ |R_{23}| &= \kappa
\\&\le |\partial (A\cap F)|\\
&=  |R_{12}|+|R_{22}|+ |R_{21}|,
\end{align*}
and hence $$|R_{23}|\le
|R_{12}|.$$

Clearly, \begin{align*}
|R_{13}|+|R_{12}|+|R_{33}|&\ge |R_{13}|+|R_{23}|+R_{33}\\&= |\nabla (F)|
\ge |A|\\&=|R_{11}|+|R_{12}|+R_{13}.\end{align*}
Thus $R_{33}\neq \emptyset.$

\begin{align*}  |R_{12}|+|R_{22}|+ |R_{32}| &=  \kappa
\\&\le |\partial (A\cup F)|\\
&\le  |R_{22}|+|R_{23}|+ |R_{32}|,
\end{align*}
and hence $|R_{12}|\le |R_{23}|.$ Thus, $|R_{12}= |R_{23}|,$ since the reverse inequality is proved above.

It follows that
   $$\kappa \le |\partial (A\cap F)|\le  |R_{12}|+|R_{22}|+ |R_{21}|\le  |R_{12}|\le |R_{23}|+|R_{22}|+ |R_{21}|=\kappa ,$$
showing that $A\cap F$ is a fragment, and hence $A\cap F=A.$\end{proof}

By a {\em basic} atom, we shall mean a $1$-atom containing $1.$ The existence of a basic atom follows since
$aA$ is an atom for any $a\in G.$
The above result has the following consequence.

\begin{proposition} \label{Cay}\cite{hejc}
Let  $ S$  be a finite generating subset of a group $G$ with $1\in S$.
If $S$ is faithful, then there is  unique basic atom of $S.$ Moreover this basic atom  is a subgroup.
\end{proposition}
\begin{proof} Assume that $S$ is faithful.
The uniqueness of the basic atom follows obviously by Theorem \ref{lem1977}.
Take an arbitrary $x\in H.$
The two atoms $H$ and $x^{-1}H$ are equal since they are basic atoms. It follows that $H^{-1}H=H.$
In particular, $H$ is a subgroup.\end{proof}

Generalizations and some applications of the last result may be found in \cite{halgebra,hast,hfour}.
We may define an atom of an arbitrary set $S$ as an atom of one translate of $S$ containing $1.$ 
As an exercise, one may prove that this notion is independent from the choice of a particular translate.
By a {\em normal}
subset of a group $G,$ we shall mean a subset closed by conjugation. The next result appears fist in \cite{hejc}.

\begin{corollary} \label{Caynormal}\cite{hejc}
Let  $ S$  be a finite generating normal subset of a group $G$ with $1\in S$.
The basic atom of $S$  is a normal subgroup of $G$ and is also  a basic atom of $S^{-1}$.
\end{corollary}
\begin{proof}
One may check easily that $xHx^{-1}$ is a basic atom of $S$ if $H$ is a basic atom of $S.$
Assume first that $S$ is faithful. By Proposition \ref{Cay}, $H$ is a subgroup and from the uniqueness of the basic atom
$xHx^{-1}=H,$ and thus $H$ is a normal subgroup. Now we have $H S^{-1}= S^{-1}H=(H S)^{-1}.$
Thus $|H S^{-1}|=|H|+\kappa (S)=|H|+\kappa ( S^{-1}).$ In particular, $H$ is also  a basic atom of $S^{-1}$.

It remains to Assume that $S$ is faithful. Suppose the contrary. Then clearly $G$ is finite. Let $K$ be a basic atom of $S^{-1}.$
By Lemma \ref{finiteg}, $S^{-1}$ is faithful. By the first part $K$ is a normal subgroup and also an atom of $S$.  Thus $|G|-|KS|-\kappa=
|G|-|SK|-\kappa\ge |K|.$ Thus $S$ is faithful, a contradiction.\end{proof}

\section{A Kneser type  result}

  \begin{theorem}\label{kneser}
Let $ S$ be a finite  faithful subset of a  group $G$  with $1\in S,$  $|S^{-1}S|\le 2|S|-2$ and $S^{-1}S\neq G.$
Let $H$ a basic atom of $S$ and let $C$ be a $H$-right-component of $S$ having a smallest possible cardinality.
Then for all $x,y\in S,$ with $(x,y)\not\in C\times C,$ one has $x^{-1}Hy\subset S^{-1}S.$

\end{theorem}
\begin{proof}

Our assumptions imply that  $\kappa (S)\le |S|-2.$
By Proposition \ref{Cay},  $H$ is a subgroup. Take a $H$-component $M$ of $S$ with $x\in M.$
Take a $H$-component $N$ of $S$ with $y\in N.$ Without loss of generality, we may take $M\neq C.$

Since $|S|-2\ge \kappa (S)=|HS|-|H|,$ we have $$2|H|-|M|-|N|\le 2|H|-|M|-|C|\le |HS|- |S|\le |H|-2.$$
Thus $|H|+2\le |M|+|N|.$
  By Lemma \ref{prehistorical},
  
$S^{-1}S \supset M^{-1}N=M^{-1}HN=x^{-1}Hy.$ \end{proof}

\begin{corollary}\label{kneserun}
Let $ A$ be a subset of group $G_0$  with $|A^{-1}A|\le 2|A|-2$ and $A^{-1}A\neq G,$
where  $G$ is the subgroup generated by ${A^{-1}A}.$

Then there are an element $a\in A$  and a non-null subgroup $H$ of $G$ such that one of the following holds:
\begin{itemize}
  \item  For all $(x,y)\in A^2\setminus (Ha)^2,$  $x^{-1}Hy  \subset A^{-1}A,$ and $|A^{-1}A|>|A^{-1}H|+|HA|-2|H|.$
  \item  For all $(x,y)\in A^2\setminus (aH)^2,$  $xHy^{-1}  \subset AA^{-1},$ and  $|AA^{-1}|>|AH|+|AH|-2|H|.$
\end{itemize}

\end{corollary}

By Lemma \ref{finiteg}, $S$ is faithful or 
 $S^{-1}$ is faithful, where $S=r^{-1}A,$ for some $r\in A.$ The reader may apply
 Theorem \ref{kneser} to get the corollary.

In Theorem \ref{kneser}, there is an uncertainty about one coset.
This uncertainty disappear in the abelian case, since the null coset has two expressions. Theorem \ref{kneser} imply Kneser's Theorem for $A-A$ and that the period of $A-A$ is the basic atom. The next corollary deals with the more general case of normal subsets of non necessarily abelian groups.

  \begin{corollary}\label{knesernormal}
Let $ 1\in S$ be a normal subset of a group $G$  with $|S^{-1}S|\le 2|S|-2$ and $S^{-1}S\neq G.$
Then the basic atom  $H$ of $S$ is a subgroup with $HS^{-1}S=S^{-1}S.$ 
\end{corollary}
\begin{proof}

By Corollary \ref{Caynormal}, $H$ is a normal subgroup. Take a $H$-component $C$  with a minimal cardinality.

Since $H$ is a proper subgroup and since $1\in S,$
there is  a $H$-component $M$ of $S$ with $C\neq M.$ Choose an element $a\in M.$ Take an arbitrary element
$x\in C.$ By Theorem \ref{kneser}, $C^{-1}C\subset a^{-1}Ha=H\subset S^{-1}S.$ By Theorem \ref{kneser},
$HS^{-1}S=S^{-1}HS=S^{-1}S.$
\end{proof}

In the abelian case, the last Corollary shows that the period of $S^{-1}S$ contains the basic atom.
As  exercise, the reader may prove equality. As a second exercise, the reader may obtain a short, using our method,
of Kneser's Theorem for $S^j,$ in the abelian case. We know no short proof based on atoms for the general form of Kneser's Theorem.

\section{Covering by a small number of cosets }
The next proposition could be useful for further investigations.

\begin{proposition} \label{lesstwothird}

Let $G$ be a group generated by a finite subset
$S$  with $1\in S$ and let $k$ be an integer such that  $|S|\ge k+1$
and $\kappa (S) <(1-\frac{1}{k})|S|.$ Then there exists a proper    subgroup $H$
such that
one of the following holds:
\begin{itemize}
  \item $ (k-1)|H|\ge |HS|\ge |S|> (k-2)|H|=\kappa (S),$
  \item  $ |HS|>(k-1)|H|\ge |SH|\ge |S|> (k-2)|H|=\kappa (S).$
\end{itemize}
 \end{proposition}

\begin{proof}

 Let
$H$  (resp. $K$) be a $1$-atom  of $S$ (resp $S^{-1}$) such that $1\in H$
 (resp. $1\in K$).

Assume first that $H$ is a subgroup
and put $|HS|=u|H|.$ By the definitions, we have
$$ \kappa (S)=|HS|-|H|=(1-\frac{1}u)u|H|\ge (1-\frac{1}u)|S|.$$
Thus $u\le k-1.$  Clearly $|S|-1=|\{1\}S|-1 \ge \kappa (S)=u|H|.$

Assume now that $H$ is not a subgroup.

By Proposition \ref{Cay},  $G$ is finite, $K$ is a subgroup and $|H| > |K|.$
By Lemma \ref{finiteg},
 $\kappa (S)= \kappa (S^{-1}).$
Now we have $|SK|-|K|=\kappa (S^{-1})=\kappa (S).$
 We must have $|KS|-|K|>\kappa (S),$ otherwise $K$ would be a $1$--atom of $S,$ a contradiction.
 The proof follows now by the first case.\end{proof}

\begin{corollary}\cite{hcras}\label{olson}
     Let  $S$  be a finite  subset of a
group $G$ with $1\in S$. Then $\kappa (S)\geq \frac{|S|}{2}.$
\end{corollary}

Corollary \ref{olson} is essentially equivalent to a result obtained independently by Olson in \cite{olsonjnt}.
The result proved in \cite{hcras}  deals with vertex-transitive graphs.
Restricted to Cayley graphs, this result   reduces to Corollary \ref{olson}.

One may find in \cite{hfour} the description of the subsets $S$ with $\kappa (S)= \frac{|S|}{2}.$
We shall use the description of the subsets $S$ with $\kappa (S)< \frac{2|S|}{3},$ given by Proposition \ref{lesstwothird}.


\section{Periodicity}

 Improving the golden ratio given in \cite{t2}, we obtain the  following result:

  \begin{theorem}\label{periodic}
Let $ A$ be a subset of group $G_0$ and let $G$ be the subgroup generated by
${A^{-1}A}.$  If $|A^{-1}A|<\min (|G|, \frac{5|A|}3),$
then there are a normal subgroup
 $K$ of  $G$ and a subgroup  $H$ with $K\subset H\subset A^{-1}A $ and  $2|K|\ge |H|$ such that
$$A^{-1}AK=KA^{-1}A=A^{-1}A \  \text{and}\   6|K|\ge|A^{-1}A|=3|H|.$$

\end{theorem}

\begin{proof}

Let us first observe that Theorem \ref{periodic} holds for a fixed subset $A$ if and only if it holds for
some translate of $A$. Take an arbitrary $a\in G$ and suppose that the result holds for a set $A.$ Since $(aA)^{-1}aA=A^{-1}A,$ Theorem \ref{periodic} holds for $A$
if and only if it holds for $aA.$  Suppose that
Theorem \ref{periodic} holds for $A.$ By Lemma \ref{Obs0}, $a^{-1}Ka$ is a normal subgroup of $\subgp{(Aa)^{-1}Aa}.$ We have
$$(Aa)^{-1}Aa=a^{-1}(A^{-1}A)a=a^{-1}A^{-1}AKa=((Aa)^{-1}Aa)(a^{-1}Ka).$$
Similarly, $(Aa)^{-1}Aa=(a^{-1}Ka)(Aa)^{-1}Aa).$
The validity of the result where  $Aa$ replaces $A$ is now obvious.

 Thus we may replace $A$ by any translate of $A.$

Put $S=r^{-1}A,$ where $r\in A.$ Since $S\subset A^{-1}A,$ we have $\subgp{S}\subset G.$ The other inclusion
 follows since $SS^{-1}=A^{-1}A.$
Notice that $1\in S$ and that $S$ generates $G.$

Observe that   $ |S^{-1}S|< \min (|G|,\frac{5|S|}3),$
and hence  $\kappa _1(S)<  \frac{2|S|}3.$ By  Proposition \ref{lesstwothird}, there is a proper subgroup $H$ such that we are in
one of the following cases:

{\bf Case} 1: $2|H|\ge |HS|\ge  |S|>|H|.$

Up to replacing $
S$ by another translate of $A,$ we may assume that
$1\in S_1$ and $|S_1|\ge |S_2|,$ where $S_1=S\cap H$ and $S_2=S\setminus H.$
Choose an element  $v\in S_2$  and put $K=(v^{-1}Hv)\cap H.$

Clearly  $2|S_1|\ge |S_1|+ |S_2|=|S|>|H|.$    By Lemma \ref{prehistorical},
$S_1^{-1}S_1=H$ and $S_1^{-1}S_2=Hv.$ In particular,

 \begin{equation}\label{mixedcosold}
 S^{-1}S=S_1^{-1}S_1\cup S_1^{-1}S_2\cup S_2^{-1}S_1\cup S_2^{-1}S_2 = H\cup v^{-1}H\cup Hv \cup S_2^{-1}S_2.
  \end{equation}

  We start by proving that
$$ 2|S_2|> |H|.$$

Assume first that $K= H,$ and hence $vH=Hv.$
  Necessarily $v^{-1}H\neq vH,$ otherwise $|S|>|H|=\frac{|G|}2,$ a contradiction.
  By Lemma \ref{Obs} and (\ref{mixedcosold}), $$ \frac{5|H|}3+\frac{5|S_2|}3\ge\frac{5|S|}3>|S^{-1}S|\ge 3|H|= \frac{5|H|}3+\frac{4|H|}3,$$ and hence $|S_2|\ge \frac{4|H|}5.$ Assume now that $K\neq H.$
  By (\ref{mixedcosold}), $$ \frac{5|H|}3+\frac{5|S_2|}3\ge \frac{5|S|}3>|S^{-1}S|\ge |Ha\setminus (a^{-1}H)|+2|H|\ge  2|H|+\frac{|H|}2,$$ and hence $2|S_2|> |H|.$

  By Lemma \ref{prehistorical}, $S_2^{-1}S_2=v^{-1}Hv.$   By (\ref{mixedcosold}), we have

  \begin{equation}\label{mixedcos}
 S^{-1}S= H\cup v^{-1}H\cup Hv \cup v^{-1}Hv.
  \end{equation}
Let us now prove that
\begin{equation}\label{2k} 2|K|\ge |H|.
  \end{equation}
We may assume that $K\neq H.$ By Lemma \ref{Obs} and  (\ref{mixedcos}), we have  $$\frac{5|S|}3>|S^{-1}S|\ge
| v^{-1}H\cup Hv|+ |H\cup  v^{-1}Hv|\ge 4|H|-2|K|\ge 3|H|+|K|(\frac{|H|}{|K|}-2),$$
and hence $|H|\le 2|K|.$

Let us prove that $K$ is a normal subgroup. If $K=H,$ then $vH=Hv.$
Since $G=\subgp{H\cup \{v\}},$ $H$ is a normal subgroup. Suppose that $K\neq H$.
We must have  $Hv\cap v^{-1}H\neq \emptyset,$  otherwise we
 have using (\ref{mixedcos}), $$|S^{-1}S|\ge  3|H|+|v^{-1}Hv\setminus H|\ge \frac{7|H|}2>\frac{5|S|}3,$$ a contradiction.
By    Lemma \ref{Obs}, $K$ is a normal subgroup of $G.$




Now we have
  using (\ref{mixedcos}),
$$KS^{-1}S=S^{-1}SK= (H\cup v^{-1}H\cup Hv \cup v^{-1}Hv)K=S^{-1}S.$$

Since $2|S_2|>|H|,$ we have $|S_2K|=|S_1K|=2|K|.$ In particular,
$|S^{-1}K|=4|K|=2|H|.$ Now we have
$ |S^{-1}S|=|KS^{-1}S|\ge |KS^{-1}|+\kappa _1(S)=2|H|+|H|=3|H|.$
We must have $ |S^{-1}S|=3|H|,$ otherwise $ |S^{-1}S|\ge 3|H|+|K|\ge \frac{7|H|}2\ge \frac{7|S|}4,$
a contradiction.

{\bf Case} 2:  $|HS|>2|H|\ge |SH|\ge  |S|>|H|.$

Up to replacing $
S$ by  $s^{-1}A,$ for some $s\in S\setminus H,$ we may assume that
$1\in S_1$ and $|S_1|\ge |S|- |S_1|,$ where $S_1=S\cap H.$ Take  a decomposition
$$S\setminus S_1=S_2\cup \dots\cup S_u,$$ where  $S_i$ is the intersection
of $S\setminus S_1$ with some right $H$-coset. We shall assume, without loss of generality, that  $|S_2|\ge \cdots\ge |S_u|.$
Since $|HS|>2|H|,$ we have  $u\ge 3.$

Since $\frac{5|S|}3>|S^{-1}S|\ge |{S_1}^{-1}S|\ge u|S_1|\ge \frac{u|S|}2,$ we have necessarily  $u=3.$

Notice that  $$|S_1|\ge  |S|- |S_1|=|S_2|+|S_3|.$$

Take $a\in S_2$ and $b\in S_3.$ Put $K=(a^{-1}Ha)\cap H.$
Notice also that $K\neq H$ and that $a \notin N(H).$  By  Lemma \ref{Obs}, we have

 \begin{equation} \label{SI}
 S_1\subset aH\cap Ha=aK \ \text{and}\ S_2\subset bK.
   \end{equation}

In particular, $|S_i|\le |K|,$ for every $i\ge 2.$

We have $$|S_1|+|S_2|\ge |H|+1,$$
otherwise we have since $|S_2|+|S_3|\le |S_1|,$

 \begin{align}|S_1^{-1}S|&= |S_1^{-1}S_1|+|S_1^{-1}S_2|+ |S_1^{-1}S_3|\nonumber\\
 &\ge |H|+2|S_1|\ge |S_1|+|S_2|+2|S_1| \nonumber
 \\&\ge\frac{5|S_1|}3+\frac{4(|S_2|+|S_3|)}3+|S_2|\ge \frac{5|S|}3,\nonumber
 \end{align}
a contradiction.
We have  \begin{equation} |S_1|+|S_3|\ge |H|+1,\label{s13}\end{equation}
otherwise
 \begin{align}|S_1^{-1}S|&= |S_1^{-1}S_1|+|S_1^{-1}S_2|+ |S_1^{-1}S_3|\nonumber
 \\&\ge 2|H|+|S_1|\ge  |S_1|+|S_3|+2|S_2|+|S_1| \nonumber\\ &=
\frac{5|S_1|}3+\frac{5|S_2|}3+\frac{|S_1|}3+\frac{|S_2|}3+|S_3|>\frac{5|S|}3,\nonumber
 \end{align}
 a contradiction.

By Lemma \ref{prehistorical}, $|S_1^{-1}S_2|=|S_1^{-1}S_3|=|H|.$
In particular,  \begin{equation}\label{mixedcoset}
S^{-1}S\supset H\cup Ha \cup Hb \cup a^{-1}H\cup {b}^{-1}H.\end{equation}

We have  \begin{equation}\label{2k2}
Ha\cap a^{-1}H\neq \emptyset, \end{equation}
 otherwise
we have using Lemma \ref{Obs} and (\ref{mixedcoset}),

 \begin{align}|S^{-1}S|&\ge |Ha \setminus b^{-1}H|+3|H|\nonumber\\
  &\ge\frac{|H|}2+ 3| H| \nonumber\ge \frac{7|H|}2\ge \frac{7|S|}4> \frac{5|S|}3,\nonumber
 \end{align}
a contradiction.

We have  $2|K|= |H|,$
otherwise we have  by Lemma \ref{Obs}, $|Ha \setminus (a^{-1}H\cup {b}^{-1}H)|\ge \frac{|H|}3,$ and hence
 \begin{align}|S^{-1}S|&\ge \frac{|H|}3+ |(a^{-1}H)\cup ( {b}^{-1}H)\cup H| \nonumber\\&\ge \frac{10|H|}3\ge \frac{5|S|}3,
 \nonumber\end{align}
a contradiction.

Thus $K$ is a normal subgroup of $H.$ By (\ref{2k2}) and Lemma \ref{Obs}, $K$ is a normal subgroup of $G.$

By (\ref{SI}),  $S_2 ^{-1}S_3$ is contained in some $K$-coset. Therefore, we must have $S_2 ^{-1}S_3\subset H\cup Ha  \cup Hb,$
otherwise $$|S^{-1}S|\ge 3|H|+|S_2|\ge 3|S_1|+|S_2|\ge \frac{5|S_1|}3+4\frac{|S_2|+|S_3|}3+|S_2|\ge \frac{5|S|}3,$$ a contradiction.
By the symmetry of $S ^{-1}S,$ we have  $S_3 ^{-1}S_2\subset H\cup Ha  \cup Hb.$
By (\ref{SI}), we have  $S_i^{-1}S_i\subset K\subset H,$ for all $i\ge 2.$

By (\ref{mixedcoset}), we have
$S^{-1}S=H\cup Ha  \cup H b.$ In particular,
$KS^{-1}S=S^{-1}SK=S^{-1}S.$
\end{proof}

{\bf Acknowledgement}
The author would like to thank  Professors Ben Green and Terence Tao, for calling his attention to Freiman's work and Husbands dissertation.

\end{document}